# Sombor indices of Silicon Carbide Graphs

Masoud Ghods[1*], Zahra Rostami[2]

1 Department of Mathematics, Semnan University, Semnan 35131-19111, Iran, mghods@semnan.ac.ir

2 Department of Mathematics, Semnan University, Semnan 35131-19111, Iran, zahrarostami.98@semnan.ac.ir

*Correspondence: mghods@semnan.ac.ir; Tel.: (09122310288)

**Abstract**

Using inequalities is a good way of studying topological indices. Chemical graph theory is one of the nontrivial applications of graph theory. In this paper, we examine and calculate another degree-based topological index for silicon-carbide structures. The Sombor index, which is a degree-based index, was introduced by Gutman in late 2020. In this paper, we first calculate and compare the Sombor index, then the decreasing Sombor index, and finally the mean Sombor index for a group of silicon-carbide structures.

## 1. Introduction

Chemical Graph theory has been playing an imperative role in mathematical chemistry. In recent years this has been the subject of much attention. For example, in [13], the first and second inverse indices of Zagreb, the first and second inverse hyper and Zagreb indices, the inverse atomic bond index and the inverse geometric-arc index for $TUC_4[m,n]$ are calculated. In [12], two families $C_n(1,4)$ and $C_n(1,5)$ of circular networks are considered and, three distance-based topological indices are calculated, namely the Wiener index, Hyper-Wiener index and a topological index. Molecular Schultz these are networks. In the literature, Kwun YC, et.al. computed some indices for the Silicon-Carbon graph. Let $G = (V_G, E_G)$ be a simple connected graph with a vertex set $V_G$ and an edge set $E_G$. The numbers of edges and vertices of $G$ are called the size $m$ and the order $n$ of $G$, respectively. The concept of silicon carbide was introduced by an American scientist in 1891. We consider a family of silicon carbide graphs, and calculate necessary computations in this article. This family includes eight silicon carbide graphs with the names and notation as: $SiC_3 - I[p,q]$, $SiC_{3\_}II[p,q]$, $SiC_{3\_}III[p,q]$, $Si_2C_3 - I[p,q]$, $Si_2C_3 - II[p,q]$, $SiC_4 - I[p,q]$, and, $SiC_4 - II[p,q]$, for arbitrary $p, q \geq 1$.

The Sombor index was first introduced by Gutman in late 2020, and since then many researchers have studied and researched it and have obtained numerous relationships based on it. To read more about the work done, you can refer to [1-11]. Sombor index is defined as[1]

$$So(G) = \sum_{uv \in E_G} \sqrt{d_u^2 + d_v^2},$$

That the $d_u$ denotes the degree of vertex $u$. We also have for the reduced Sombor index and an average Sombor index, respectively:

$$So_{red}(G) = \sum_{uv \in E_G} \sqrt{(d_u - 1)^2 + (d_v - 1)^2},$$

$$So_{avr}(G) = \sum_{uv \in E_G} \sqrt{(d_u - \frac{2m}{n})^2 + (d_v - \frac{2m}{n})^2}.$$

We consider a family of silicon carbide graphs, and calculate necessary computations in this article. This family includes eight silicon carbide graphs with the names and notation as: $SiC_3 - I[p,q]$, $SiC_{3\_} II[p,q]$, $SiC_{3\_} III[p,q]$, $Si_2C_3 - I[p,q]$, $Si_2C_3 - II[p,q]$, $SiC_4 - I[p,q]$, and, $SiC_4 - II[p,q]$, for arbitrary $p, q \geq 1$.

1. **Mane result**

In this section, in Figures 1–8, carbon atoms are shown as red, and silicon atoms $Si$ are shown as blue.

1. In the Silicon-Carbon ($SiC_{3\_} I[p,q]$), the total number of vertices are $8pq$, and the total number of edges are $12pq - 2p - 3q$. Since, we have:

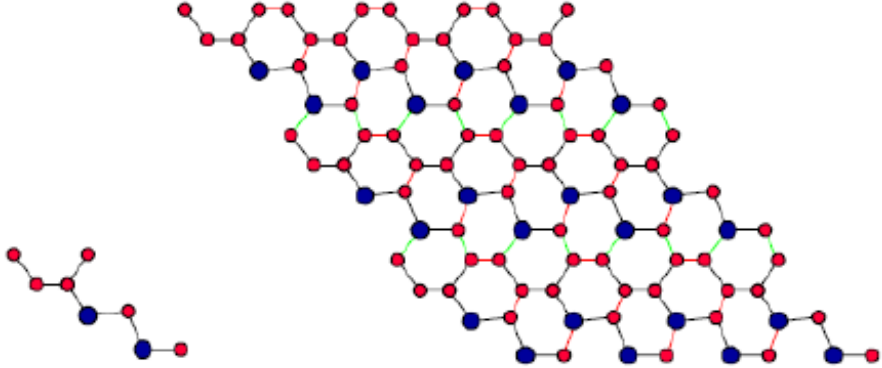

**Figure 1.** The Unit cell of $SiC_3 - I[p,q]$, The Sheet of $SiC_3 - I[p,q]$ for $p = 4$, $q = 3$.

$$E_{\{1,2\}}(SiC_3 - I[p,q]) = \{e = uv \in E(SiC_3 - I[p,q]) | d_u = 1, d_v = 2\},$$
$$E_{\{1,3\}}(SiC_3 - I[p,q]) = \{e = uv \in E(SiC_3 - I[p,q]) | d_u = 1, d_v = 3\},$$
$$E_{\{2,2\}}(SiC_3 - I[p,q]) = \{e = uv \in E(SiC_3 - I[p,q]) | d_u = 2, d_v = 2\},$$
$$E_{\{2,3\}}(SiC_3 - I[p,q]) = \{e = uv \in E(SiC_3 - I[p,q]) | d_u = 2, d_v = 3\},$$

Additionally, it has

$$E_{\{3,3\}}(SiC_3 - I[p,q]) = \{e = uv \in E(SiC_3 - I[p,q]) | d_u = 3, d_v = 3\}.$$

Now, we have

$$|E_{\{1,2\}}(SiC_3 - I[p,q])| = 2,$$
$$|E_{\{1,3\}}(SiC_3 - I[p,q])| = 1,$$
$$|E_{\{2,2\}}(SiC_3 - I[p,q])| = \begin{cases} 3q - 1 & p = 1, q \geq 1 \\ 2p + 2q - 3 & p > 1, q \geq 1 \end{cases},$$
$$|E_{\{2,3\}}(SiC_3 - I[p,q])| = \begin{cases} 6q - 4 & p = 1, q \geq 1 \\ 4p + 8q - 8 & p > 1, q \geq 1 \end{cases},$$

And

$$|E_{\{3,3\}}(SiC_3 - I[p,q])| = \begin{cases} 12pq - 2p - 12q + 2 & p = 1, q \geq 1 \\ 12pq - 8q - 13p + 8 & p > 1, q \geq 1 \end{cases}.$$

**Theorem 1** consider $G = SiC_3 - I[p,q]$ be the Silicon Carbide. Then we have

$$So(SiC_3 - I[p,q])$$
$$= \begin{cases} 36\sqrt{2}pq + (6\sqrt{13} - 30\sqrt{2})q - 6\sqrt{2}p + (4\sqrt{2} + 2\sqrt{5} + \sqrt{10} - 4\sqrt{13}) & p = 1, q \geq 1 \\ 36\sqrt{2}pq + (8\sqrt{13} - 20\sqrt{2})q + (4\sqrt{13} - 35\sqrt{2})p + (18\sqrt{2} + 2\sqrt{5} + \sqrt{10} - 8\sqrt{13}) & p > 1, q \geq 1 \end{cases}.$$

***Proof.*** Let $G$ has parameters $p$ and $q$. we can see that the total number of vertices are $8pq$, and the total number of edges is $15pq - 2p - 3q$. By sombor index, we have:

For $p = 1, q \geq 1$:

$$So(SiC_3 - I[p,q]) = \sum_{uv \in E_G} \sqrt{d_u^2 + d_v^2}$$
$$= 2\sqrt{1^2 + 2^2} + \sqrt{1^2 + 3^2} + (3q - 1)\sqrt{2^2 + 2^2} + (6q - 4)\sqrt{2^2 + 3^2}$$
$$+ (12pq - 2p - 12q + 2)\sqrt{3^2 + 3^2}$$
$$= 2\sqrt{5} + \sqrt{10} + 2\sqrt{2}(3q - 1) + (6q - 4)\sqrt{13} + (12pq - 2p - 12q + 2)3\sqrt{2}$$
$$= 36\sqrt{2}pq + (6\sqrt{13} - 30\sqrt{2})q - 6\sqrt{2}p + (4\sqrt{2} + 2\sqrt{5} + \sqrt{10} - 4\sqrt{13}).$$

For $p > 1, q \geq 1$:

$$So(SiC_3 - I[p,q]) = \sum_{uv \in E_G} \sqrt{d_u^2 + d_v^2}$$
$$= 2\sqrt{1^2 + 2^2} + \sqrt{1^2 + 3^2} + (2p + 2q - 3)\sqrt{2^2 + 2^2} + (4p + 8q - 8)\sqrt{2^2 + 3^2}$$
$$+ (12pq - 8q - 13p + 8)\sqrt{3^2 + 3^2}$$
$$= 2\sqrt{5} + \sqrt{10} + 2\sqrt{2}(2p + 2q - 3) + (4p + 8q - 8)\sqrt{13}$$
$$+ (12pq - 8q - 13p + 8)3\sqrt{2}$$
$$= 36\sqrt{2}pq + (8\sqrt{13} - 20\sqrt{2})q + (4\sqrt{13} - 35\sqrt{2})p + (18\sqrt{2} + 2\sqrt{5} + \sqrt{10} - 8\sqrt{13}).$$

□

**Theorem 2** consider $G = SiC_3 - I[p,q]$ be the Silicon Carbide. Then we have

$$So_{red}(SiC_3 - I[p,q])$$
$$= \begin{cases} 24\sqrt{2}pq + (6\sqrt{5} - 21\sqrt{2})q - 4\sqrt{2}p + (4 + 3\sqrt{2} - 4\sqrt{5}) & p = 1, q \geq 1 \\ 24\sqrt{2}pq + (8\sqrt{5} - 14\sqrt{2})q + (4\sqrt{5} - 24\sqrt{2})p + (4 - 8\sqrt{5} + 13\sqrt{2}) & p > 1, q \geq 1 \end{cases}.$$

***Proof.*** For $p = 1, q \geq 1$:

$$So(SiC_3 - I[p,q]) = \sum_{uv \in E_G} \sqrt{(d_u - 1)^2 + (d_v - 1)^2}$$
$$= 2\sqrt{(1-1)^2 + (2-1)^2} + \sqrt{(1-1)^2 + (3-1)^2}$$
$$+ (3q - 1)\sqrt{(2-1)^2 + (2-1)^2} + (6q - 4)\sqrt{(2-1)^2 + (3-1)^2}$$
$$+ (12pq - 2p - 12q + 2)\sqrt{(3-1)^2 + (3-1)^2}$$
$$= 2\sqrt{1} + 2 + \sqrt{2}(3q - 1) + (6q - 4)\sqrt{5} + (12pq - 2p - 12q + 2)2\sqrt{2}$$
$$= 24\sqrt{2}pq + (6\sqrt{5} - 21\sqrt{2})q - 4\sqrt{2}p + (4 + 3\sqrt{2} - 4\sqrt{5}).$$

For $p > 1, q \geq 1$:

$$So(SiC_3 - I[p,q]) = \sum_{uv \in E_G} \sqrt{(d_u - 1)^2 + (d_v - 1)^2}$$
$$= 2\sqrt{(1-1)^2 + (2-1)^2} + \sqrt{(1-1)^2 + (3-1)^2}$$
$$+ (2p + 2q - 3)\sqrt{(2-1)^2 + (2-1)^2} + (4p + 8q - 8)\sqrt{(2-1)^2 + (3-1)^2}$$
$$+ (12pq - 8q - 13p + 8)\sqrt{(3-1)^2 + (3-1)^2}$$
$$= 2\sqrt{1} + \sqrt{4} + \sqrt{2}(2p + 2q - 3) + (4p + 8q - 8)\sqrt{5}$$
$$+ (12pq - 8q - 13p + 8)2\sqrt{2}$$
$$= 24\sqrt{2}pq + (8\sqrt{5} - 14\sqrt{2})q + (4\sqrt{5} - 24\sqrt{2})p + (4 - 8\sqrt{5} + 13\sqrt{2}).$$

□

2. In the Silicon-Carbon ($SiC_{3\_}II[p,q]$), the total number of vertices are $8pq$, and the total number of edges are $12pq - 2p - 2q$. Since, we have:

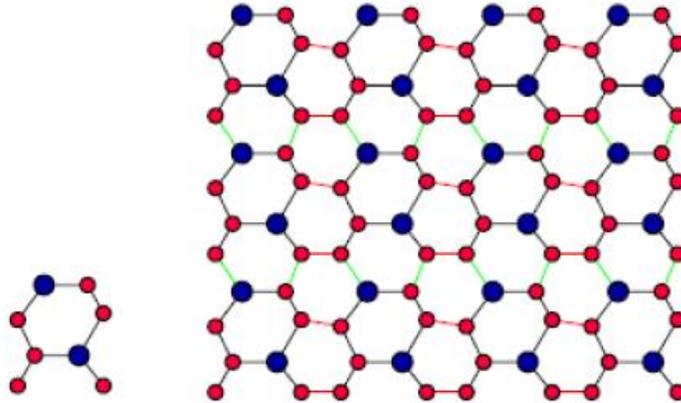

**Figure 1.** The Unit cell of $SiC_3 - II[p,q]$, The Sheet of $SiC_3 - II[p,q]$ for $p = 4, q = 3$.

$$E_{\{1,3\}}(SiC_3 - II[p,q]) = \{e = uv \in E(SiC_3 - II[p,q]) | d_u = 1, d_v = 3\},$$

$$E_{\{2,2\}}(SiC_3 - II[p,q]) = \{e = uv \in E(SiC_3 - II[p,q]) | d_u = 2, d_v = 2\},$$
$$E_{\{2,3\}}(SiC_3 - II[p,q]) = \{e = uv \in E(SiC_3 - II[p,q]) | d_u = 2, d_v = 3\},$$

Additionally, it has
$$E_{\{3,3\}}(SiC_3 - I[p,q]) = \{e = uv \in E(SiC_3 - II[p,q]) | d_u = 3, d_v = 3\}.$$

Now, we have
$$|E_{\{1,3\}}(SiC_3 - II[p,q])| = 2,$$
$$|E_{\{2,2\}}(SiC_3 - II[p,q])| = 2p + 1,$$
$$|E_{\{2,3\}}(SiC_3 - II[p,q])| = 4p + 8q - 10,$$

and
$$|E_{\{3,3\}}(SiC_3 - II[p,q])| = 12pq - 8p - 10q + 7.$$

**Theorem 3** consider $G = SiC_3 - II[p,q]$ be the Silicon Carbide. Then we have

$$\begin{aligned}So(SiC_3 - II[p,q]) \\ = 36\sqrt{2}pq + (8\sqrt{13} - 30\sqrt{2})q + (4\sqrt{13} - 20\sqrt{2})p \\ + (2\sqrt{10} - 10\sqrt{13} + 23\sqrt{2}).\end{aligned}$$

*Proof.* Let $G$ has parameters $p$ and $q$. we can see that the total number of vertices are $8pq$, and the total number of edges is $12pq - 2p - 2q$. Then it can be easily proved by placing the degree of vertices of each edge and the number of edges in relation to the sombor index.

□

**Theorem 4** consider $G = SiC_3 - II[p,q]$ be the Silicon Carbide. Then we have

$$\begin{aligned}So_{red}(SiC_3 - II[p,q]) \\ = 36\sqrt{2}pq + (8\sqrt{5} - 20\sqrt{2})q + (4\sqrt{5} - 14\sqrt{2})p + (4 - 10\sqrt{5} + 15\sqrt{2}).\end{aligned}$$

*Proof.* Let $G$ has parameters $p$ and $q$. we can see that the total number of vertices are $8pq$, and the total number of edges is $12pq - 2p - 2q$. Then it can be easily proved by placing the degree of vertices of each edge and the number of edges in relation to the reduced sombor index.

□

3. In the Silicon-Carbon ($SiC_{3\_}III[p,q]$), the total number of vertices are $8pq$, and the total number of edges are $12pq - 2p - 3q$. Since, we have:

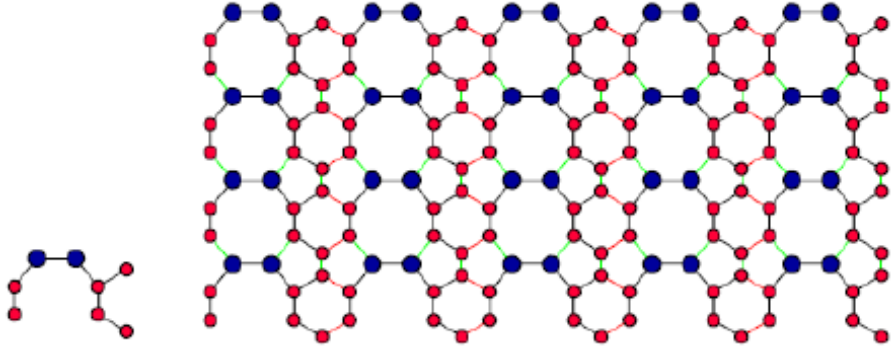

**Figure 3.** The Unit cell of $SiC_3 - III[p,q]$, The Sheet of $SiC_3 - III[p,q]$ for $p = 5$, $q = 4$.

$$E_{\{1,2\}}(SiC_3 - III[p,q]) = \{e = uv \in E(SiC_3 - III[p,q]) | d_u = 1, d_v = 2\},$$
$$E_{\{1,3\}}(SiC_3 - III[p,q]) = \{e = uv \in E(SiC_3 - III[p,q]) | d_u = 1, d_v = 3\},$$
$$E_{\{2,2\}}(SiC_3 - III[p,q]) = \{e = uv \in E(SiC_3 - III[p,q]) | d_u = 2, d_v = 2\},$$
$$E_{\{2,3\}}(SiC_3 - III[p,q]) = \{e = uv \in E(SiC_3 - III[p,q]) | d_u = 2, d_v = 3\},$$

Additionally, it has
$$E_{\{3,3\}}(SiC_3 - III[p,q]) = \{e = uv \in E(SiC_3 - III[p,q]) | d_u = 3, d_v = 3\}.$$

Now, we have
$$|E_{\{1,2\}}(SiC_3 - III[p,q])| = 1,$$
$$|E_{\{1,3\}}(SiC_3 - III[p,q])| = 2,$$
$$|E_{\{2,2\}}(SiC_3 - III[p,q])| = 3p + 2q - 3,$$
$$|E_{\{2,3\}}(SiC_3 - III[p,q])| = 6p + 4q - 8,$$

and

$$|E_{\{3,3\}}(SiC_3 - III[p,q])| = 12pq - 12p - 8q + 8.$$

**Theorem 5** consider $G = SiC_3 - III[p,q]$ be the Silicon Carbide. Then we have

$$So(SiC_3 - III[p,q])$$
$$= 36\sqrt{2}pq + (4\sqrt{13} - 22\sqrt{2})q + (6\sqrt{13} - 30\sqrt{2})p$$
$$+ (2\sqrt{10} - 8\sqrt{13} + 18\sqrt{2} + \sqrt{5}).$$

*Proof.* Let $G$ has parameters $p$ and $q$. we can see that the total number of vertices are $8pq$, and the total number of edges is $12pq - 2p - 3q$. Then it can be easily proved by placing the degree of vertices of each edge and the number of edges in relation to the sombor index.

□

**Theorem 6** consider $G = SiC_3 - III[p,q]$ be the Silicon Carbide. Then we have

$$So_{red}(SiC_3 - III[p,q])$$
$$= 24\sqrt{2}pq + (4\sqrt{5} - 14\sqrt{2})q + (6\sqrt{5} - 21\sqrt{2})p + (5 - 8\sqrt{5} + 13\sqrt{2}).$$

**Proof.** Let $G$ has parameters $p$ and $q$. we can see that the total number of vertices are $8pq$, and the total number of edges is $12pq - 2p - 3q$. Then it can be easily proved by placing the degree of vertices of each edge and the number of edges in relation to the reduced sombor index.

□

4. In the Silicon-Carbon ($Si_2C_3 - I[p,q]$), the total number of vertices are $10pq$, and the total number of edges are $15pq - 2p - 3q$. Since, we have:

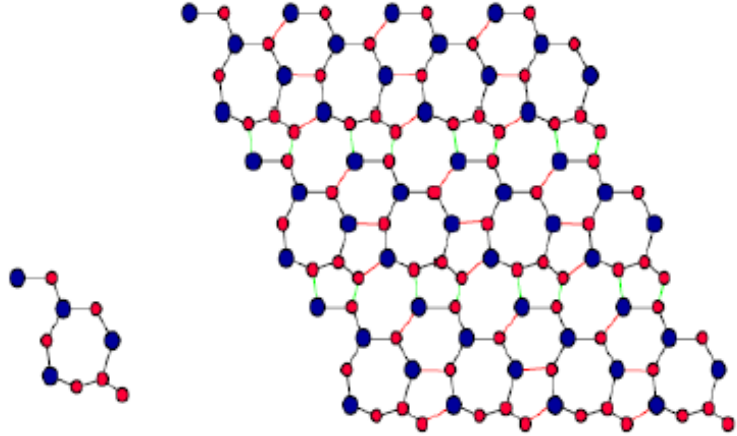

**Figure 4.** The Unit cell of $Si_2C_3 - I[p,q]$, The Sheet of $Si_2C_3 - I[p,q]$ for $p = 4$, $q = 3$.

$$E_{\{1,2\}}(Si_2C_3 - I[p,q]) = \{e = uv \in E(Si_2C_3 - I[p,q]) | d_u = 1, d_v = 2\},$$
$$E_{\{1,3\}}(Si_2C_3 - I[p,q]) = \{e = uv \in E(Si_2C_3 - I[p,q]) | d_u = 1, d_v = 3\},$$
$$E_{\{2,2\}}(Si_2C_3 - I[p,q]) = \{e = uv \in E(Si_2C_3 - I[p,q]) | d_u = 2, d_v = 2\},$$
$$E_{\{2,3\}}(Si_2C_3 - I[p,q]) = \{e = uv \in E(Si_2C_3 - I[p,q]) | d_u = 2, d_v = 3\},$$

Additionally, it has
$$E_{\{3,3\}}(Si_2C_3 - I[p,q]) = \{e = uv \in E(Si_2C_3 - I[p,q]) | d_u = 3, d_v = 3\}.$$

Now, we have
$$|E_{\{1,2\}}(Si_2C_3 - I[p,q])| = 1$$
$$|E_{\{1,3\}}(Si_2C_3 - I[p,q])| = 1,$$
$$|E_{\{2,2\}}(Si_2C_3 - I[p,q])| = p + 2q,$$
$$|E_{\{2,3\}}(Si_2C_3 - I[p,q])| = 6p - 1 + 8(q - 1),$$

And
$$|E_{\{3,3\}}(Si_2C_3 - I[p,q])| = 15pq - 9p - 13q + 7.$$

**Theorem 7** consider $G = Si_2C_3 - II[p,q]$ be the Silicon Carbide. Then we have

$$So(Si_2C_3 - I[p,q])$$
$$= 45\sqrt{2}pq + (8\sqrt{13} - 35\sqrt{2})q + (6\sqrt{13} - 25\sqrt{2})p$$
$$+ (21\sqrt{2} + \sqrt{5} + \sqrt{10} - 9\sqrt{13}).$$

*Proof.* Let $G$ has parameters $p$ and $q$. we can see that the total number of vertices are $10pq$, and the total number of edges is $15pq - 2p - 3q$. Then it can be easily proved by placing the degree of vertices of each edge and the number of edges in relation to the sombor index.

□

**Theorem 8** consider $G = Si_2C_3 - II[p,q]$ be the Silicon Carbide. Then we have

$$So_{red}(Si_2C_3 - I[p,q])$$
$$= 30\sqrt{2}pq + (8\sqrt{5} - 16\sqrt{2})q + (6\sqrt{5} - 17\sqrt{2})p + (3 + 14\sqrt{2} - 9\sqrt{5}).$$

*Proof.* Let $G$ has parameters $p$ and $q$. we can see that the total number of vertices are $10pq$, and the total number of edges is $15pq - 2p - 3q$. Then it can be easily proved by placing the degree of vertices of each edge and the number of edges in relation to the reduced sombor index.

□

5. In the Silicon-Carbon $(Si_2C_3 - II[p,q])$, the total number of vertices are $10pq$, and the total number of edges are $15pq - 2p - 3q$. Since, we have:

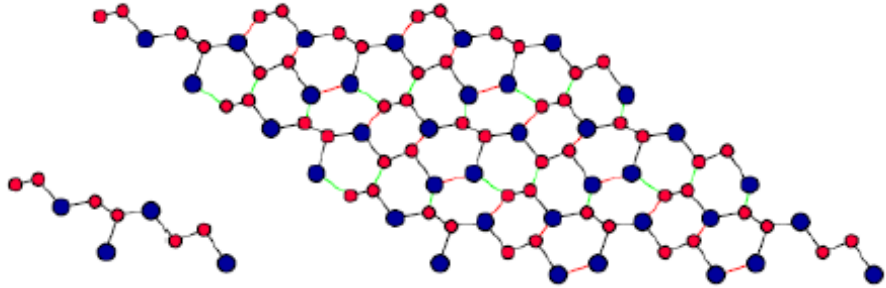

**Figure5.** The Unit cell of $Si_2C_3 - II[p,q]$, The Sheet of $Si_2C_3 - II[p,q]$ for $p = 3$, $q = 3$.

$E_{\{1,2\}}(Si_2C_3 - II[p,q]) = \{e = uv \in E(Si_2C_3 - II[p,q]) | d_u = 1, d_v = 2\}$,
$E_{\{1,3\}}(Si_2C_3 - II[p,q]) = \{e = uv \in E(Si_2C_3 - II[p,q]) | d_u = 1, d_v = 3\}$,
$E_{\{2,2\}}(Si_2C_3 - II[p,q]) = \{e = uv \in E(Si_2C_3 - II[p,q]) | d_u = 2, d_v = 2\}$,
$E_{\{2,3\}}(Si_2C_3 - II[p,q]) = \{e = uv \in E(Si_2C_3 - II[p,q]) | d_u = 2, d_v = 3\}$,

Additionally, it has

$E_{\{3,3\}}(Si_2C_3 - II[p,q]) = \{e = uv \in E(Si_2C_3 - II[p,q]) | d_u = 3, d_v = 3\}$.

Now, we have

$$|E_{\{1,2\}}(Si_2C_3 - II[p,q])| = 2$$
$$|E_{\{1,3\}}(Si_2C_3 - II[p,q])| = 1,$$

$$|E_{\{2,2\}}(Si_2C_3 - II[p,q])| = 2p + 2q,$$

$$|E_{\{2,3\}}(Si_2C_3 - II[p,q])| = 8p + 8q - 14,$$

And

$$|E_{\{3,3\}}(Si_2C_3 - II[p,q])| = 15pq - 13p - 13q + 11.$$

**Theorem 9** consider $G = Si_2C_3 - II[p,q]$ be the Silicon Carbide. Then we have
$So(Si_2C_3 - II[p,q])$
$$= 45\sqrt{2}pq + (8\sqrt{13} - 35\sqrt{2})q + (8\sqrt{13} - 35\sqrt{2})p + (35\sqrt{2} + 2\sqrt{5} + \sqrt{10} - 14\sqrt{13})$$

***Proof.*** Let $G$ has parameters $p$ and $q$. we can see that the total number of vertices are $10pq$, and the total number of edges is $15pq - 2p - 3q$. Then it can be easily proved by placing the degree of vertices of each edge and the number of edges in relation to the sombor index.
□

**Theorem 10** consider $G = Si_2C_3 - II[p,q]$ be the Silicon Carbide. Then we have
$So_{red}(Si_2C_3 - II[p,q])$
$$= 30\sqrt{2}pq + (8\sqrt{5} - 24\sqrt{2})q + (8\sqrt{5} - 24\sqrt{2})p + (3 + 22\sqrt{2} - 14\sqrt{5}).$$

***Proof.*** Let $G$ has parameters $p$ and $q$. we can see that the total number of vertices are $10pq$, and the total number of edges is $15pq - 2p - 3q$. Then it can be easily proved by placing the degree of vertices of each edge and the number of edges in relation to the reduced sombor index.
□

6. In the Silicon-Carbon ($Si_2C_3 - III[p,q]$), the total number of vertices are $10pq$, and the total number of edges are $15pq - 2p - 3q$. Since, we have:

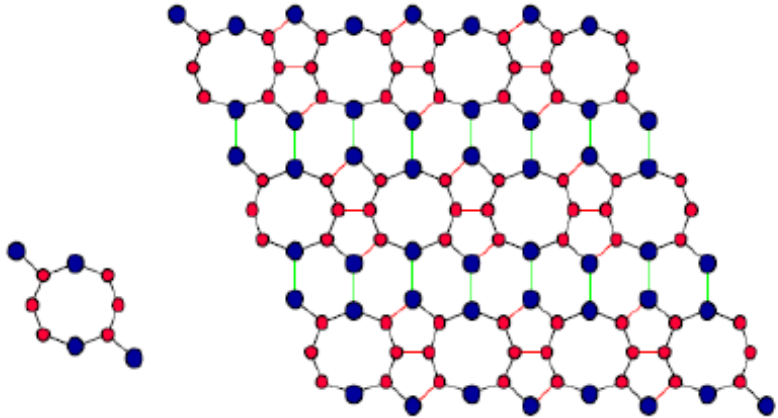

**Figure 6.** The Unit cell of $Si_2C_3 - III[p,q]$, The Sheet of $Si_2C_3 - III[p,q]$ for $p = 4$, $q = 3$.

$$E_{\{1,3\}}(Si_2C_3 - III[p,q]) = \{e = uv \in E(Si_2C_3 - III[p,q]) | d_u = 1, d_v = 3\},$$
$$E_{\{2,2\}}(Si_2C_3 - III[p,q]) = \{e = uv \in E(Si_2C_3 - III[p,q]) | d_u = 2, d_v = 2\},$$

$$E_{\{2,3\}}(Si_2C_3 - III[p,q]) = \{e = uv \in E(Si_2C_3 - III[p,q])|d_u = 2, d_v = 3\},$$
Additionally, it has
$$E_{\{3,3\}}(Si_2C_3 - III[p,q]) = \{e = uv \in E(Si_2C_3 - III[p,q])|d_u = 3, d_v = 3\}.$$
Now, we have
$$|E_{\{1,3\}}(Si_2C_3 - III[p,q])| = 2,$$
$$|E_{\{2,2\}}(Si_2C_3 - III[p,q])| = 2 + 2q,$$
$$|E_{\{2,3\}}(Si_2C_3 - III[p,q])| = 8p + 8q - 12,$$
And
$$|E_{\{3,3\}}(Si_2C_3 - II[p,q])| = 15pq - 10p - 13q + 8.$$

**Theorem 11** consider $G = Si_2C_3 - III[p,q]$ be the Silicon Carbide. Then we have
$$So(Si_2C_3 - III[p,q])$$
$$= 45\sqrt{2}pq + (8\sqrt{13} - 35\sqrt{2})q + (8\sqrt{13} - 30\sqrt{2})p$$
$$+ (28\sqrt{2} + 2\sqrt{10} - 12\sqrt{13}).$$

**Proof.** Let $G$ has parameters $p$ and $q$. we can see that the total number of vertices are $10pq$, and the total number of edges is $15pq - 2p - 3q$. Then it can be easily proved by placing the degree of vertices of each edge and the number of edges in relation to the sombor index.
□

**Theorem 12** consider $G = Si_2C_3 - II[p,q]$ be the Silicon Carbide. Then we have
$$So_{red}(Si_2C_3 - III[p,q])$$
$$= 30\sqrt{2}pq + (8\sqrt{5} - 24\sqrt{2})q + (8\sqrt{5} - 20\sqrt{2})p + (4 + 18\sqrt{2} - 12\sqrt{5}).$$

**Proof.** Let $G$ has parameters $p$ and $q$. we can see that the total number of vertices are $10pq$, and the total number of edges is $15pq - 2p - 3q$. Then it can be easily proved by placing the degree of vertices of each edge and the number of edges in relation to the reduced sombor index.
□

7. In the Silicon-Carbon ($SiC_4 - I[p,q]$), the total number of vertices are $10pq$, and the total number of edges are $15pq - 4p - 2q + 1$. Since, we have:

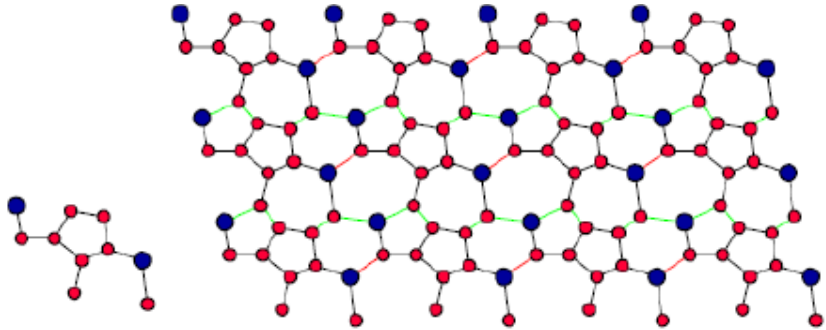

**Figure 7.** The Unit cell of $SiC_4 - I[p,q]$, The Sheet of $SiC_4 - I[p,q]$ for $p = 4$, $q = 3$.

$$E_{\{1,2\}}(SiC_4 - I[p,q]) = \{e = uv \in E(SiC_4 - I[p,q])|d_u = 1, d_v = 2\},$$
$$E_{\{1,3\}}(SiC_4 - I[p,q]) = \{e = uv \in E(SiC_4 - I[p,q])|d_u = 1, d_v = 3\},$$
$$E_{\{2,2\}}(SiC_4 - I[p,q]) = \{e = uv \in E(SiC_4 - I[p,q])|d_u = 2, d_v = 2\},$$
$$E_{\{2,3\}}(SiC_4 - I[p,q]) = \{e = uv \in E(SiC_4 - I[p,q])|d_u = 2, d_v = 3\},$$

Additionally, it has
$$E_{\{3,3\}}(SiC_4 - I[p,q]) = \{e = uv \in E(SiC_4 - I[p,q])|d_u = 3, d_v = 3\}.$$

Now, we have
$$|E_{\{1,2\}}(SiC_4 - I[p,q])| = 2,$$
$$|E_{\{1,3\}}(SiC_4 - I[p,q])| = 3p - 2,$$
$$|E_{\{2,2\}}(SiC_4 - I[p,q])| = p + 2q - 2,$$
$$|E_{\{2,3\}}(SiC_4 - I[p,q])| = 2p + 4q - 2,$$

And
$$|E_{\{3,3\}}(SiC_4 - I[p,q])| = 14pq - 10p - 8q + 5.$$

**Theorem 13** consider $G = SiC_4 - II[p,q]$ be the Silicon Carbide. Then we have

$$So(Si_2C_3 - I[p,q])$$
$$= 42\sqrt{2}pq + (4\sqrt{13} - 20\sqrt{2})q + (2\sqrt{13} + 3\sqrt{10} - \sqrt{2})p$$
$$+ (11\sqrt{2} + 2\sqrt{5} - 2\sqrt{10} - 2\sqrt{13}).$$

*Proof.* Let $G$ has parameters $p$ and $q$. we can see that the total number of vertices are $10pq$, and the total number of edges is $15pq - 4p - 2q + 1$. Then it can be easily proved by placing the degree of vertices of each edge and the number of edges in relation to the sombor index.

□

**Theorem 14** consider $G = SiC_4 - I[p,q]$ be the Silicon Carbide. Then we have

$$So_{red}(SiC_4 - II[p,q])$$
$$= 28\sqrt{2}pq + (4\sqrt{5} - 14\sqrt{2})q + (6 + 2\sqrt{5} - 14\sqrt{2})p + (8\sqrt{2} - 2\sqrt{5} - 2).$$

*Proof.* Let $G$ has parameters $p$ and $q$. we can see that the total number of vertices are $10pq$, and the total number of edges is $15pq - 4p - 2q + 1$. Then it can be easily proved by placing the degree of vertices of each edge and the number of edges in relation to the reduced sombor index.

□

8. In the Silicon-Carbon ($SiC_4 - II[p,q]$), the total number of vertices are $10pq$, and the total number of edges are $15pq - 4p - 2q$. Since, we have:

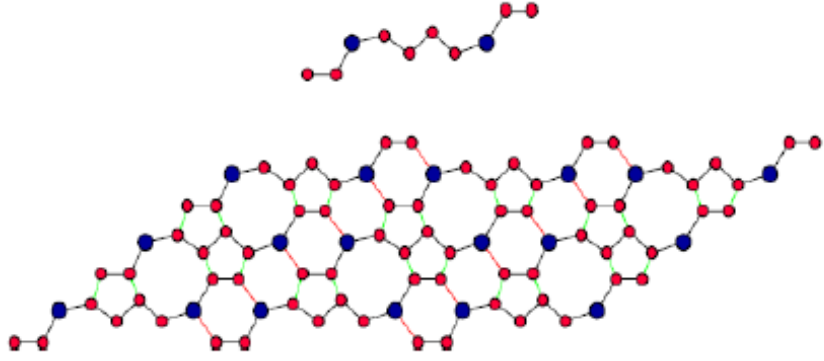

**Figure 8.** The Unit cell of $SiC_4 - II[p,q]$, The Sheet of $SiC_4 - II[p,q]$ for $p = 3$, $q = 3$.

$$E_{\{1,2\}}(SiC_4 - II[p,q]) = \{e = uv \in E(SiC_4 - II[p,q]) | d_u = 1, d_v = 2\},$$
$$E_{\{2,2\}}(SiC_4 - II[p,q]) = \{e = uv \in E(SiC_4 - II[p,q]) | d_u = 2, d_v = 2\},$$
$$E_{\{2,3\}}(SiC_4 - II[p,q]) = \{e = uv \in E(SiC_4 - II[p,q]) | d_u = 2, d_v = 3\},$$

Additionally, it has

$$E_{\{3,3\}}(SiC_4 - II[p,q]) = \{e = uv \in E(SiC_4 - II[p,q]) | d_u = 3, d_v = 3\}.$$

Now, we have

$$|E_{\{1,2\}}(SiC_4 - II[p,q])| = 2,$$
$$|E_{\{2,2\}}(SiC_4 - II[p,q])| = \begin{cases} 5q + 2 & q = 1, p \geq 1 \\ 2p + 2 & q > 1, p \geq 1 \end{cases},$$
$$|E_{\{2,3\}}(SiC_4 - II[p,q])| = \begin{cases} 6p - 6 & q = 1, p \geq 1 \\ 12p + 8q - 14 & q > 1, p \geq 1 \end{cases}$$

And

$$|E_{\{3,3\}}(SiC_4 - II[p,q])| = \begin{cases} 15pq - 15p - 2q & q = 1, p \geq 1 \\ 12pq - 18q - 10p & q > 1, p \geq 1 \end{cases}.$$

**Theorem 15** consider $G = SiC_4 - II[p,q]$ be the Silicon Carbide. Then we have

$$So(SiC_4 - II[p,q])$$
$$= \begin{cases} 45\sqrt{2}pq - 6\sqrt{2}q + (6\sqrt{13} - 35\sqrt{2}p) + (4\sqrt{2} + 2\sqrt{5} - 6\sqrt{13}) & p = 1, q \geq 1 \\ 45\sqrt{2}pq + (8\sqrt{13} - 30\sqrt{2})q + (12\sqrt{13} - 50\sqrt{2})p + (4\sqrt{2} + 2\sqrt{5} - 14\sqrt{13}) & p > 1, q \geq 1 \end{cases}.$$

*Proof.* Let $G$ has parameters $p$ and $q$. we can see that the total number of vertices are $10pq$, and the total number of edges is $15pq - 4p - 2q$. Then it can be easily proved by placing the degree of vertices of each edge and the number of edges in relation to the sombor index.

□

**Theorem 16** consider $G = SiC_4 - II[p,q]$ be the Silicon Carbide. Then we have

$$So_{red}(SiC_4 - II[p,q])$$
$$= \begin{cases} 30\sqrt{2}pq + (-4\sqrt{2})q + (6\sqrt{5} - 25\sqrt{2})p + (2 + 2\sqrt{2} - 6\sqrt{5}) & p = 1, q \geq 1 \\ 30\sqrt{2}pq + (8\sqrt{5} - 20\sqrt{2})q + (12\sqrt{5} - 34\sqrt{2})p + (2 - 14\sqrt{5} + 2\sqrt{2}) & p > 1, q \geq 1 \end{cases}.$$

**Proof.** Let $G$ has parameters $p$ and $q$. we can see that the total number of vertices are $10pq$, and the total number of edges is $15pq - 4p - 2q$. Then it can be easily proved by placing the degree of vertices of each edge and the number of edges in relation to the reduced sombor index.

□

**Theorem 17** consider $G$ be the Silicon Carbide. Then we have

$$So_{Avg}(G) = \sqrt{2}\left(|E_1|\sqrt{(A - 3/2)^2 + 1/4} + |E_2|\sqrt{(A - 1)^2 + 1} + |E_3|(A - 1)\right.$$
$$\left. + |E_4|\sqrt{(A - 1/2)^2 + 1/4} + |E_5|A\right).$$

Such that

| $G$ | $\|E_1\|$ | $\|E_2\|$ | $\|E_3\|$ | | $\|E_4\|$ | | $\|E_5\|$ | | $A$ |
|---|---|---|---|---|---|---|---|---|---|
| $SiC_3 - I[p,q]$ | 2 | 1 | $2q - 1$ $\quad p=1, q \geq 1$<br>$2p + 2q - 3 \quad p > 1, q \geq 1$ | | $6q - 4 \quad p=1, q \geq 1$<br>$4p + 8q - 8 \quad p > 1, q \geq 1$ | | $12pq - 2p - 12q + 2 \quad p=1, q \geq 1$<br>$12pq - 8q - 13p + 8 \quad p > 1, q \geq 1$ | | $\dfrac{2p + 3q}{4pq}$ |
| $SiC_3 - II[p,q]$ | 0 | 2 | $2p + 1$ | | $4p + 8q - 10$ | | $12pq - 8p - 10q + 7$ | | $\dfrac{p + q}{2pq}$ |
| $SiC_3 - III[p,q]$ | 1 | 2 | $3p + 2q - 3$ | | $6p + 4q - 8$ | | $12pq - 12p - 8q + 8$ | | $\dfrac{3p + 2q}{4pq}$ |
| $Si_2C_3 - I[p,q]$ | 1 | 1 | $p + 2q$ | | $6p + 8q - 9$ | | $15pq - 9p - 13q + 7$ | | $\dfrac{2p + 3q}{5pq}$ |
| $Si_2C_3 - II[p,q]$ | 2 | 1 | $2p + 2q$ | | $8p + 8q - 14$ | | $15pq - 13p - 13q + 11$ | | $\dfrac{3p + 3q}{5pq}$ |
| $Si_2C_3 - III[p,q]$ | 0 | 2 | $2 + 2q$ | | $8p + 8q - 12$ | | $15pq - 10p - 13q + 8$ | | $\dfrac{2p + 3q}{5pq}$ |
| $SiC_4 - I[p,q]$ | 2 | $3p - 2$ | $p + 2q - 2$ | | $2p + 4q - 2$ | | $14pq - 10p - 8q + 5$ | | $\dfrac{4p + q - 1}{5pq}$ |
| $SiC_4 - II[p,q]$ | 2 | 0 | $5q + 2 \quad q=1, p \geq 1$<br>$2p + 2 \quad q > 1, p \geq 1$ | | $6p - 6 \quad q=1, p \geq 1$<br>$12p + 8q - 14 \quad q > 1, p \geq 1$ | | $15pq - 15p - 2q \quad q=1, p \geq 1$<br>$12pq - 18q - 10p \quad q > 1, p \geq 1$ | | $\dfrac{4p + 2q}{5pq}$ |

**Proof:** consider $G$ be $SiC_3 - I[p, q]$. Then

$$|E_1| = |E_{\{1,2\}}(SiC_3 - I[p,q])| = 2,$$
$$|E_2| = |E_{\{1,3\}}(SiC_3 - I[p,q])| = 1,$$
$$|E_3| = |E_{\{2,2\}}(SiC_3 - I[p,q])| = \begin{cases} 3q - 1 & p = 1, q \geq 1 \\ 2p + 2q - 3 & p > 1, q \geq 1' \end{cases}$$
$$|E_4| = |E_{\{2,3\}}(SiC_3 - I[p,q])| = \begin{cases} 6q - 4 & p = 1, q \geq 1 \\ 4p + 8q - 8 & p > 1, q \geq 1' \end{cases}$$

And

$$|E_5| = |E_{\{3,3\}}(SiC_3 - I[p,q])| = \begin{cases} 12pq - 2p - 12q + 2 & p = 1, q \geq 1 \\ 12pq - 8q - 13p + 8 & p > 1, q \geq 1. \end{cases}$$

Now, using the definition of the average Sombor index, we have:

$$So_{avr}(G) = \sum_{uv \in E_G} \sqrt{(d_u - \frac{2m}{n})^2 + (d_v - \frac{2m}{n})^2}.$$

That $\frac{2m}{n} = \frac{2(12pq-2p-3q)}{8pq} = \frac{24pq-4p-6q}{8pq} = 3 - (\frac{2p+3q}{4pq})$. Now by setting the equation we get:

$$So_{avr}(SiC_3 - I[p,q]) = \sum_{uv \in E_G} \sqrt{(d_u - (3 - (\frac{2p+3q}{4pq})))^2 + (d_v - (3 - (\frac{2p+3q}{4pq})))^2}$$

**Case 1:** For $p = 1, q \geq 1$:

$$So(SiC_3 - I[p,q]) = \sum_{uv \in E_G} \sqrt{\left(d_u - \left(3 - \left(\frac{2p+3q}{4pq}\right)\right)\right)^2 + \left(d_v - \left(3 - \left(\frac{2p+3q}{4pq}\right)\right)\right)^2}$$

$$= 2\sqrt{\left(1 - \left(3 - \left(\frac{2p+3q}{4pq}\right)\right)\right)^2 + \left(2 - \left(3 - \left(\frac{2p+3q}{4pq}\right)\right)\right)^2}$$

$$+ \sqrt{\left(1 - \left(3 - \left(\frac{2p+3q}{4pq}\right)\right)\right)^2 + \left(3 - \left(3 - \left(\frac{2p+3q}{4pq}\right)\right)\right)^2}$$

$$+ (3q - 1)\sqrt{\left(2 - \left(3 - \left(\frac{2p+3q}{4pq}\right)\right)\right)^2 + \left(2 - \left(3 - \left(\frac{2p+3q}{4pq}\right)\right)\right)^2}$$

$$+ (6q - 4)\sqrt{\left(2 - \left(3 - \left(\frac{2p+3q}{4pq}\right)\right)\right)^2 + \left(3 - \left(3 - \left(\frac{2p+3q}{4pq}\right)\right)\right)^2}$$

$$+ (12pq - 2p - 12q + 2)\sqrt{\left(3 - \left(3 - \left(\frac{2p+3q}{4pq}\right)\right)\right)^2 + \left(3 - \left(3 - \left(\frac{2p+3q}{4pq}\right)\right)\right)^2}$$

$$= 2\sqrt{\left(-2+\left(\frac{2p+3q}{4pq}\right)\right)^2 + \left(-1+\left(\frac{2p+3q}{4pq}\right)\right)^2} + \sqrt{\left(-2+\left(\frac{2p+3q}{4pq}\right)\right)^2 + \left(0+\left(\frac{2p+3q}{4pq}\right)\right)^2}$$

$$+ (3q-1)\sqrt{\left(-1+\left(\frac{2p+3q}{4pq}\right)\right)^2 + \left(-1+\left(\frac{2p+3q}{4pq}\right)\right)^2}$$

$$+ (6q-4)\sqrt{\left(-1+\left(\frac{2p+3q}{4pq}\right)\right)^2 + \left(0+\left(\frac{2p+3q}{4pq}\right)\right)^2}$$

$$+ (12pq - 2p - 12q + 2)\sqrt{\left(0+\left(\frac{2p+3q}{4pq}\right)\right)^2 + \left(0+\left(\frac{2p+3q}{4pq}\right)\right)^2}$$

Using the binomial square union, we rewrite the above relation:

$So(SiC_3 - I[p,q])$

$$= 2\sqrt{4 - 4\left(\frac{2p+3q}{4pq}\right) + \left(\frac{2p+3q}{4pq}\right)^2 + 1 - 2\left(\frac{2p+3q}{4pq}\right) + \left(\frac{2p+3q}{4pq}\right)^2}$$

$$+ \sqrt{4 - 4\left(\frac{2p+3q}{4pq}\right) + \left(\frac{2p+3q}{4pq}\right)^2 + \left(\frac{2p+3q}{4pq}\right)^2}$$

$$+ (3q-1)\sqrt{1 - 2\left(\frac{2p+3q}{4pq}\right) + \left(\frac{2p+3q}{4pq}\right)^2 + 1 - 2\left(\frac{2p+3q}{4pq}\right) + \left(\frac{2p+3q}{4pq}\right)^2}$$

$$+ (6q-4)\sqrt{1 - 2\left(\frac{2p+3q}{4pq}\right) + \left(\frac{2p+3q}{4pq}\right)^2 + \left(\frac{2p+3q}{4pq}\right)^2}$$

$$+ (12pq - 2p - 12q + 2)\sqrt{\left(\frac{2p+3q}{4pq}\right)^2 + \left(\frac{2p+3q}{4pq}\right)^2}$$

$$= 2\sqrt{2\left(\frac{2p+3q}{4pq}\right)^2 - 6\left(\frac{2p+3q}{4pq}\right) + 5} + \sqrt{2\left(\frac{2p+3q}{4pq}\right)^2 - 4\left(\frac{2p+3q}{4pq}\right) + 4}$$

$$+ (3q-1)\sqrt{2\left(\frac{2p+3q}{4pq}\right)^2 - 4\left(\frac{2p+3q}{4pq}\right) + 2}$$

$$+ (6q-4)\sqrt{2\left(\frac{2p+3q}{4pq}\right)^2 - 2\left(\frac{2p+3q}{4pq}\right) + 1}$$

$$+ (12pq - 2p - 12q + 2)\sqrt{2\left(\frac{2p+3q}{4pq}\right)^2}$$

$$= \sqrt{2}\left(2\sqrt{\left(\frac{2p+3q}{4pq}-\frac{3}{2}\right)^2+\frac{1}{4}}+\sqrt{\left(\frac{2p+3q}{4pq}-1\right)^2+1}+(3q-1)\left(\frac{2p+3q}{4pq}-1\right)\right.$$

$$\left.+(6q-4)\sqrt{\left(\frac{2p+3q}{4pq}-\frac{1}{2}\right)^2+\frac{1}{4}}+(12pq-2p-12q+2)\left(\frac{2p+3q}{4pq}\right)\right).$$

**Case 2:** For $p > 1, q \geq 1$:

The same case 1 is easily provable:

$$So(SiC_3-I[p,q]) = \sum_{uv\in E_G}\sqrt{\left(d_u-\left(3-\left(\frac{2p+3q}{4pq}\right)\right)\right)^2+\left(d_v-\left(3-\left(\frac{2p+3q}{4pq}\right)\right)\right)^2}$$

$$= 2\sqrt{\left(1-\left(3-\left(\frac{2p+3q}{4pq}\right)\right)\right)^2+\left(2-\left(3-\left(\frac{2p+3q}{4pq}\right)\right)\right)^2}$$

$$+\sqrt{\left(1-\left(3-\left(\frac{2p+3q}{4pq}\right)\right)\right)^2+\left(3-\left(3-\left(\frac{2p+3q}{4pq}\right)\right)\right)^2}$$

$$+(2p+2q-3)\sqrt{\left(2-\left(3-\left(\frac{2p+3q}{4pq}\right)\right)\right)^2+\left(2-\left(3-\left(\frac{2p+3q}{4pq}\right)\right)\right)^2}$$

$$+(4p+8q-8)\sqrt{\left(2-\left(3-\left(\frac{2p+3q}{4pq}\right)\right)\right)^2+\left(3-\left(3-\left(\frac{2p+3q}{4pq}\right)\right)\right)^2}$$

$$+(12pq-8q-13p+8)\sqrt{\left(3-\left(3-\left(\frac{2p+3q}{4pq}\right)\right)\right)^2+\left(3-\left(3-\left(\frac{2p+3q}{4pq}\right)\right)\right)^2}$$

$$= \sqrt{2}\left(2\sqrt{\left(\frac{2p+3q}{4pq}-\frac{3}{2}\right)^2+\frac{1}{4}}+\sqrt{\left(\frac{2p+3q}{4pq}-1\right)^2+1}\right.$$

$$+(2p+2q-3)\left(\frac{2p+3q}{4pq}-1\right)+(4p+8q-8)\sqrt{\left(\frac{2p+3q}{4pq}-\frac{1}{2}\right)^2+\frac{1}{4}}$$

$$\left.+(12pq-8q-13p+8)\left(\frac{2p+3q}{4pq}\right)\right).$$

The proof is complete.

□

**Theorem 18** consider $G$ be the Silicon Carbide. Then we have

$$So(G) \geq \frac{1}{2}M_1(G), \qquad So(G) \geq \frac{1}{3}M_2(G), \qquad So(G) \geq 2ISI(G)$$

Where $M_1(G), M_2(G)$, and $ISI(G)$ are denoted as the first Zagreb index, second Zagreb index, and inverse sum indeg index, respectively.

*Proof:* For g we prove that the other cases are similarly provable

$$M_1(Si_2C_3 - I[p,q]) = \sum_{uv \in E_G} d_u + d_v$$
$$= 1(1+2) + 1(1+3) + (p+2q)(2+2) + (6p+8q-9)(2+3)$$
$$+ (15pq - 9p - 13q + 7)(3+3) = 90pq - 20p - 30q + 4. \qquad (1)$$

$$M_2(Si_2C_3 - I[p,q]) = \sum_{uv \in E_G} d_u \times d_v$$
$$= 1(1 \times 2) + 1(1 \times 3) + (p+2q)(2 \times 2) + (6p+8q-9)(2 \times 3)$$
$$+ (15pq - 9p - 13q + 7)(3 \times 3) = 135pq - 41p - 61q + 14. \qquad (2)$$

$$ISI(Si_2C_3 - I[p,q]) = \sum_{uv \in E_G} d_u \times d_v \times (d_u + d_v)^{-1} \cong 22.5pq - 5.3p - 7.9q + 1.11. \qquad (3)$$

$$So(Si_2C_3 - I[p,q])$$
$$= 45\sqrt{2}pq + \left(8\sqrt{13} - 35\sqrt{2}\right)q + \left(6\sqrt{13} - 25\sqrt{2}\right)p$$
$$+ \left(21\sqrt{2} + \sqrt{5} + \sqrt{10} - 9\sqrt{13}\right) \cong 42.45pq - 10.63p - 16.1q + 38.45. \qquad (4)$$

By comparing relationships 1-4, the results are easily obtained.

□

2. **Conclusion**

This article is a continuation of our previous studies [14]. In this paper, we obtain Sombor indices for silicon-carbide. We also examined the general case for the Average Sombor index. In the following, we will discuss the boundaries for the Sombor indices and other indices obtained for silicon carbide. As you can see in this article, one of these boundaries is presented in the last theorem.

**Author contributions:** the authors made equal contributions in the article.
**Acknowledgments:** we thanks the reviewers for their suggestions, which have helped us to improve the quality of this article.
**Conflicts of interest:** the authors declare no conflict of interest.